\documentclass[12pt,reqno]{amsart}

\usepackage{amsmath,amsfonts,amsthm,amssymb,amsxtra}
\usepackage{bbm} % to get \1
\usepackage{hyperref}	% to get hyperlinks in equations and references
\newtheorem{theorem}{Theorem}%[section]

\newtheorem{lemma}[theorem]{Lemma}
\newtheorem{corollary}[theorem]{Corollary}

\theoremstyle{definition}

\newtheorem{definition}[theorem]{Definition}

\theoremstyle{remark}

\newtheorem{remark}[theorem]{Remark}

\numberwithin{equation}{section}
\numberwithin{theorem}{section}

\newcommand{\C}{\mathbb{C}}

\renewcommand{\epsilon}{\varepsilon}

\renewcommand{\phi}{\varphi}
\newcommand{\R}{\mathbb{R}}

\newcommand{\Z}{\mathbb{Z}}

\def\key{\vspace{.5em}
{\textit{Keywords}:\,\relax%
}}

\def\MSC{\vspace{.5em}
{2020 \textit{Mathematics Subject Classification}:\,\relax%
}}

\begin{document}

\title[Hardy and Rellich inequality on lattices]{Hardy and Rellich inequality on lattices}

\author{By Shubham Gupta\\}
\address[Shubham Gupta]{Department of Mathematics, Imperial College London, 180 Queen’s Gate, London,
SW7 2AZ, United Kingdom.}
\email{s.gupta19@imperial.ac.uk}

\begin{abstract}
In this paper, we study the asymptotic behaviour of the sharp constant in discrete Hardy and Rellich inequality on the lattice $\Z^d$ as $d \rightarrow \infty$. In the process, we proved some Hardy-type inequalities for the operators $\Delta^m$ and $\nabla(\Delta^m)$ for non-negative integers $m$ on a $d$ dimensional torus. It turns out that the sharp constant in discrete Hardy and Rellich inequality grows as $d$ and $d^2$ respectively as $ d \rightarrow \infty$.
\end{abstract}
\maketitle

\key{Discrete Hardy inequality, Discrete Rellich inequality, Hardy inequality on torus, Rellich inequality on torus, Expansion of squares, Sharp constant.}\\\\ 
\MSC{39B62, 26D10, 35A23.}

\section{Introduction and Main Results}
In this paper, we are interested in the asymptotic behaviour of sharp constants (largest constant for which an inequality holds true) in the following inequalities as $d \rightarrow \infty$:
\begin{equation}\label{1.1}
    \sum_{n \in \Z^d} |Du(n)|^2 \geq C_H(d) \sum_{n \in \Z^d} \frac{|u(n)|^2}{|n|^2},
\end{equation}
and 
\begin{equation}\label{1.2}
    \sum_{n \in \Z^d} |\Delta u(n)|^2 \geq C_R(d) \sum_{n \in \Z^d} \frac{|u(n)|^2}{|n|^4},
\end{equation}
where $Du(n) := (D_1u(n), D_2u(n),.., D_d u(n))$, 
\begin{align*}
    D_ju(n) := u(n)-u(n-e_j), \hspace{9pt} \Delta u(n) := \sum_{j=1}^d 2u(n)-u(n-e_j)-u(n+e_j),
\end{align*}
and $e_j$ is the $j^{th}$ canonical basis of $\R^d$.\\

Inequalities \eqref{1.1} and \eqref{1.2} are discrete analogues of Hardy and Rellich inequalities respectively on $\R^d$ (see \cite{GL} and references therein), and will be referred to as discrete Hardy and Rellich inequalities. The discrete Hardy inequality \eqref{1.1} has been considered in the past in works \cite{KL, RS1, RS2}and more generally for graphs in \cite{KPP1}. To our best knowledge \cite{KPP3} is the only paper where \eqref{1.2} has been studied in the past in the context of graphs. In particular, they prove inequality \eqref{1.2} with a weight that grows like $|n|^{-4}$ as $|n| \rightarrow \infty$ for $d \geq 5$. \\

Although Hardy and Rellich type inequalities in the continuum are very well studied and the literature on them is enormous(for instance see book by Balinsky et al. \cite{BEL}), very little is known about their discrete counterparts. In fact, the sharp constant in \eqref{1.1} and \eqref{1.2} is only known when $d=1$. Computation of sharp constant in \eqref{1.1} goes back to G.H. Hardy \cite{HLP}(see \cite{KPP2} for a more recent proof and improvement of \eqref{1.1} for $d=1$). For \eqref{1.2}, sharp constant was recently obtained in $d=1$ \cite{GKS}. In fact, authors in \cite{GKS} proved an improvement of sharp Rellich inequality \eqref{1.2} by adding lower order terms in the RHS of \eqref{1.2}. \\

In this work, we compute the behaviour of sharp constants in these inequalities for large dimensions. In fact, we will also find the asymptotic behaviour of sharp constants in the higher order versions of inequalities \eqref{1.1} and \eqref{1.2}. In the following theorems $C_c(\Z^d)$ denotes the space of finitely supported functions on $\Z^d$. The main results of the paper can be summarized as follows:

\begin{theorem}\label{thm1.1}
Let $k \geq 0$ and let $u \in C_c(\Z^d)$ with $u(0)=0$. Let $C_1(k, d)$ be the sharp constant in the following inequality:
\begin{equation}\label{1.3}
    \sum_{n \in \Z^d} |D(\Delta^k u)(n)|^2 \geq C_1(k, d) \sum_{n \in \Z^d} \frac{|u(n)|^2}{|n|^{4k+2}}.
\end{equation}
Then $C_1(k, d) \sim d^{2k+1}$ as $d \rightarrow \infty$, that is, there exists positive constants $c_1, c_2$ and $N$ independent of $d$ such that $c_1 d^{2k+1}\leq C_1(k,d) \leq c_2 d^{2k+1}$ for all $ d\geq N$. 
\end{theorem}

\begin{theorem}\label{thm1.2}
Let $k \geq 1$ and let $u \in C_c(\Z^d)$ with $u(0)=0$. Let $C_2(k, d)$ be the sharp constant in the following inequality:
\begin{equation}\label{1.4}
    \sum_{n \in \Z^d} |\Delta^k u(n)|^2 \geq C_2(k, d) \sum_{n \in \Z^d} \frac{|u(n)|^2}{|n|^{4k}}.
\end{equation}
Then $C_2(k, d) \sim d^{2k}$ as $d \rightarrow \infty$.
\end{theorem}

\begin{remark}
We would like to point out that the sharp constants in continuous analogues of \eqref{1.3} and \eqref{1.4} on $\R^d$ grows as $d^{4k+2}$ and $d^{4k}$ as $d \rightarrow \infty$ respectively, see \cite{DH}, for the computation of sharp constants.   
\end{remark}

\begin{remark}
Note that putting $k=0$ in \eqref{1.3} and $k=1$ in \eqref{1.4} give the discrete Hardy \eqref{1.1} and Rellich inequality \eqref{1.2} respectively. 
\end{remark} 

\begin{corollary}
Let $u \in C_c(\Z^d)$ with $u(0) =0$. Let $C_H(d)$ be the sharp constant in the discrete Hardy inequality
\begin{equation}\label{1.5}
    \sum_{n \in \Z^d} |Du(n)|^2 \geq C_H(d) \sum_{n \in \Z^d} \frac{|u(n)|^2}{|n|^2}.
\end{equation}
Then $C_H(d) \sim d$, as $d \rightarrow \infty$.
\end{corollary}

\begin{corollary}
Let $u \in C_c(\Z^d)$ with $u(0) =0$. Let $C_R(d)$ be the sharp constant in the discrete Rellich inequality
\begin{equation}
    \sum_{n \in \Z^d} |\Delta u(n)|^2 \geq C_R(d) \sum_{n \in \Z^d} \frac{|u(n)|^2}{|n|^4}.
\end{equation}
Then $C_R(d) \sim d^2$, as $d \rightarrow \infty$.
\end{corollary}

\begin{remark}
We would like to mention that in \cite{KPP1} authors proved \eqref{1.5} with $C_H(d)|n|^{-2}$ replaced by an optimal Hardy weight $w(n)$ (see \cite{KPP1} for the definition of optimal Hardy weight), which behaves as $w(n) = 4^{-1}(d-2)^2 |n|^{-2} +$ lower order terms, under the limit $|n| \rightarrow \infty$ . It is interesting to note that the constant in the leading term of optimal Hardy weight grows as $d^2$ as compared to linear growth of sharp constant in \eqref{1.5}.
\end{remark}

\begin{remark}\label{rem1.8}
Consider a function $u$ defined on $\Z^d$ as follows: $u(n) = 1$ when $|n|=1$ and zero otherwise. Using this function in inequalities \eqref{1.3} and \eqref{1.4} gives us the bounds $C_1(k, d) \leq 4^{2k+1}d^{2k+1} \sim d^{2k+1}$ and $C_2(k, d) \leq 4^{2k}d^{2k} \sim d^{2k}$(this will be proved in section \ref{sec:mainproofs}). We will prove inequalities \eqref{1.3}, \eqref{1.4} with explicit constants which grow as $d^{2k+1}$ and $d^{2k}$ respectively, as $d \rightarrow \infty$, thereby completing the proofs of our main results. 
\end{remark}

\begin{remark}
A natural question that arises from this work is the determination of exact constant in the leading term of the asymptotic expansion of $C_H(d)$ as $d \rightarrow \infty$, in other words what is the value of $\lim\limits_{d \rightarrow \infty} C_H(d)/d ?$ We could not answer this question and it remains open.
\end{remark}

In the next section we convert inequalities \eqref{1.3}, \eqref{1.4} into continuous Hardy type inequalities on a torus. In section \ref{sec:hardytorus}, we will prove various higher order Hardy-type inequalities on the torus, which we believe are completely new. Finally in section \ref{sec:mainproofs} we will use the results proved in sections \ref{sec:conversion} and \ref{sec:hardytorus} to prove the main results of the paper.

\section{Equivalent integral inequalities}\label{sec:conversion}
In this section, we will convert the inequalities \eqref{1.3} and \eqref{1.4} into some equivalent integral inequalities on the torus. In this paper, $Q_d := (-\pi, \pi)^d$ denotes the open square in $\R^d$. 
\begin{definition}
Let $\psi : \overline{Q_d} \rightarrow \C$ be a map. Then we say it is $2\pi$-periodic in each variable if 
\begin{align*}
    \psi(x_1,..,x_{i-1}, -\pi, x_{i+1},.., x_d) = \psi(x_1,..,x_{i-1}, \pi, x_{i+1},.., x_d),
\end{align*}
for all $1 \leq i \leq d$.
\end{definition}
\begin{lemma}\label{lem2.2}
Let $k \geq 0$ be an integer. Let $ u\in C_c(\Z^d)$ with $u(0)=0$. There exists $\psi \in C^\infty(\overline{Q_d})$, all of whose derivatives are $2\pi$-periodic in each variable and which has zero average, such that
\begin{equation}\label{2.1}
    \sum_{n \in \Z^d} \frac{|u(n)|^2}{|n|^{4k+2}} = \int_{Q_d} |\nabla(\Delta^k \psi)(x)|^2 dx,
\end{equation}
and 
\begin{equation}\label{2.2}
    \sum_{n \in \Z^d} |D(\Delta^k u)(n)|^2 = 4^{2k+1}\int_{Q_d} |\Delta^{2k+1} \psi(x)|^2 \omega(x)^{2k+1}dx,
\end{equation}
where 
\begin{equation}\label{2.3}
    \omega(x) := \sum_{j=1}^d \sin^2(x_j/2).
\end{equation}
\end{lemma}

\begin{proof}
We will first prove the result for $k=0$ and then extend the proof to general $k$. \\  

Let $u \in l^2(\Z^d)$. We define its Fourier transform $\widehat{u}$ as 
\begin{align*}
    \widehat{u}(x) := (2\pi)^{-\frac{d}{2}} \sum_{n \in \Z^d} u(n) e^{-i n \cdot x}, \hspace{9pt} x \in (-\pi, \pi)^d.
\end{align*}
Let $u_j(n) := \frac{n_j}{|n|^2} u(n) $ for $n \neq 0$ and $u_j(0) =0$. Then Parseval's identity gives us
\begin{equation}
    \sum_{n \in \Z^d}\frac{|u(n)|^2}{|n|^2} = \sum_{j=1}^d \sum_{n \in \Z
    ^d}| u_j|^2 = \int_{Q_d} \sum_{j=1}^d |\widehat{u_j}|^2 dx.
\end{equation}
Using the inversion formula for Fourier transform we get
\begin{align*}
    u(n) - u(n-e_j) = (2\pi)^{-\frac{d}{2}}\int_{Q_d} \widehat{u}(x)(1-e^{-ix_j}) e^{in \cdot x} dx.
\end{align*}
Applying Parseval's identity and summing w.r.t. to $j$, we obtain
\begin{equation}
    \begin{split}
        \sum_{n \in \Z^d} \sum_{j=1}^d |u(n)-u(n-e_j)|^2 &= 4\int_{Q_d} |\widehat{u}(x)|^2 \sum_{j=1}^d  \sin^2(x_j/2)\\
        &= 4\int_{Q_d} \Big|\sum_{j=1}^d \partial_{x_j}\widehat{u_j}(x)\Big|^2 \sum_{j=1}^d  \sin^2(x_j/2) dx,    
    \end{split}
\end{equation}
where the last identity uses $\sum_j \partial_{x_j}\widehat{u_j}(x) = -i \widehat{u}$.\\

The inversion formula for Fourier transform along with integration by parts gives us
\begin{align*}
    (2\pi)^{\frac{d}{2}}n_k u_j(n) = \int_{Q_d} n_k \widehat{u_j}(x) e^{i n \cdot x} = i \int_{Q_d}  \partial_{x_k} \widehat{u_j}(x) e^{i n \cdot x}, 
\end{align*}
which implies that
\begin{equation}
    \partial_{x_k}\widehat{u_j} = \partial_{x_j}\widehat{u_k}.
\end{equation}
This further implies that there exists a smooth function $\psi$ such that $\widehat{u_j}(x) = \partial_{x_j} \psi(x)$, whose average is zero. It is easy to see that periodicity of $\widehat{u_j}$ along with its zero average imply that $\psi$ is also $2\pi$ periodic in each variable. This proves the result for $k=0$.\\

Next we prove the result for general $k > 0$ in a similar way. Let $j_1, j_2,.., j_{2k+1}$ be integers lying between $1$ and $d$. Consider the family of functions $u_{j_1j_2..j_{2k+1}}(n) := \frac{n_{j_1}..n_{j_{2k+1}}}{|n|^{4k+2}}$ for $n \neq 0$ and $u_{j_1j_2..j_{2k+1}}(0):=0$. Then we have
\begin{equation}\label{2.7}
    \sum_{n \in \Z^d} \frac{|u(n)|^2}{|n|^{4k+2}} = \sum_{1 \leq j_1,..,j_{2k+1} \leq d} \sum_{n \in \Z^d} |u_{j_1..j_{2k+1}}(n)|^2 =  \sum_{1 \leq j_1,..,j_{2k+1} \leq d} \int_{Q_d} |\widehat{u_{j_1..j_{2k+1}}}(x)|^2 dx.
\end{equation}
Using the inversion formula for the Fourier transform, we obtain 
\begin{equation}\label{2.8}
    \widehat{D_j u}(x) = (1-e^{-ix_j})\widehat{u}(x) \hspace{5pt} \text{and} \hspace{5pt} \widehat{\Delta u}(x) = 4 \sum_{j=1}^d \sin^2(x_j/2) \widehat{u}(x).
\end{equation}
Expression \eqref{2.8} along with Parseval's identity yields
\begin{equation}\label{2.9}
    \begin{split}
        \sum_j \sum_{n \in \Z^d} |D_j \Delta^k u(n)|^2 &= \sum_j \int_{Q_d} |\widehat{D_j \Delta^k u}|^2 dx\\
        &= 4^{2k+1}\int_{Q_d} |\widehat{u}|^2 \Big(\sum_j \sin^2(x_j/2)\Big)^{2k+1}\\
        &= 4^{2k+1} \int_{Q_d} \Big|\sum_{1 \leq j_1,...,j_{2k+1}\leq d} \partial_{x_{j_{2k+1}}}..\partial_{x_{j_1}}\widehat{u_{j_1...j_{2k+1}}}\Big|^2  \omega^{2k+1}.   
    \end{split}
\end{equation}
We used $\sum_{1 \leq j_1,..,j_{2k+1}\leq d}\partial_{x_{j_{2k+1}}}..\partial_{x_{j_1}}\widehat{u_{j_1...j_{2k+1}}} = (-i)^{2k+1}\widehat{u}(x)$ in the last identity. Using the inversion formula and integration by parts, we further notice that for fixed $j_2,.., j_{2k+1}$
\begin{align*}
    \partial_{x_k} \widehat{u_{jj_2..j_{2k+1}}} = \partial_{x_j} \widehat{u_{kj_2..j_{2k+1}}}.
\end{align*}
This implies that $u_{j_1j_2..j_{2k+1}} = \partial_{x_{j_1}}\varphi_{j_2..j_{2k+1}}$ for some smooth function $\varphi_{j_2...j_{2k+1}}$ with zero average. Furthermore, along with its zero average, the periodicity of $\widehat{u_{j_1..j_{2k+1}}}$ implies that $\varphi_{j_2..j_{2k+1}}$ is $2\pi$-periodic in each variable. Since, $u_{j_1j_2..j_{2k+1}}$ is symmetric w.r.t. to $j_1$ and $j_2$ we get $\partial_{x_k}\varphi_{jj_3..j_{2k+1}} = \partial_{x_j} \varphi_{kj_3..j_{2k+1}}$. This implies that $\varphi_{j_2..j_{2k+1}}= \partial_{x_{j_2}}\xi_{j_3..j_{2k+1}}$ for some smooth $2\pi$-periodic function $\xi_{j_3...j_{2k+1}}$ whose average is zero. Using this argument iteratively, we get a smooth $2\pi$-periodic function $\psi(x)$ with zero average such that
\begin{equation}\label{2.10}
    \widehat{u_{j_1..j_{2k+1}}}(x) = \partial_{x_{j_{2k+1}}}...\partial_{x_{j_1}} \psi(x).
\end{equation}
Using \eqref{2.10} in \eqref{2.7} and \eqref{2.9}, we get expressions \eqref{2.1} and \eqref{2.2} respectively.
\end{proof}

\begin{lemma}\label{lem2.3}
Let $k \geq 0$ be an integer. Let $ u\in C_c(\Z^d)$ with $u(0)=0$. There exists $\psi \in C^\infty(\overline{Q_d})$, all of whose derivatives are $2\pi$-periodic in each variable and which has zero average, such that
\begin{equation}\label{2.11}
    \sum_{n \in \Z^d} \frac{|u(n)|^2}{|n|^{4k}} = \int_{Q_d} |\Delta^k \psi(x)|^2 dx,
\end{equation}
and 
\begin{equation}\label{2.12}
    \sum_{n \in \Z^d} |\Delta^k u(n)|^2 = 4^{2k}\int_{Q_d} |\Delta^{2k} \psi(x)|^2 \omega(x)^{2k}dx,
\end{equation}
where 
\begin{equation}\label{2.3}
    \omega(x) := \sum_{j=1}^d \sin^2(x_j/2).
\end{equation} 
\end{lemma}

\begin{proof}
The proof of this Lemma follows the proof of Lemma \ref{lem2.2} step by step but we include the proof here for the sake of completeness.\\

Let $1 \leq j_1, j_2,.., j_{2k} \leq d$. Consider the family of functions $u_{j_1j_2..j_{2k}}(n) := \frac{n_{j_1}..n_{j_{2k}}}{|n|^{4k}}$ for $n \neq 0$ and $u_{j_1j_2..j_{2k}}(0):=0$. Then we have
\begin{equation}\label{2.14}
    \sum_{n \in \Z^d} \frac{|u(n)|^2}{|n|^{4k}} = \sum_{1 \leq j_1,..,j_{2k} \leq d} \sum_{n \in \Z^d} |u_{j_1..j_{2k}}(n)|^2 =  \sum_{1 \leq j_1,..,j_{2k} \leq d} \int_{Q_d} |\widehat{u_{j_1..j_{2k}}}(x)|^2 dx.
\end{equation}
Using the inversion formula for the Fourier transform we obtain 
\begin{equation}\label{2.15}
    \widehat{\Delta u}(x) = 4 \sum_{j=1}^d \sin^2(x_j/2) \widehat{u}(x).
\end{equation}
Expressions \eqref{2.15} along with Parseval's identity yields
\begin{equation}\label{2.16}
    \begin{split}
        \sum_{n \in \Z^d} |\Delta^k u(n)|^2 = \int_{Q_d} |\widehat{\Delta^k u}|^2 dx &= 4^{2k}\int_{Q_d} |\widehat{u}|^2 \Big(\sum_j \sin^2(x_j/2)\Big)^{2k}\\
        &= 4^{2k} \int_{Q_d} \Big|\sum_{1 \leq j_1,...,j_{2k}\leq d} \partial_{x_{j_{2k}}}..\partial_{x_{j_1}}\widehat{u_{j_1...j_{2k}}}\Big|^2  \omega^{2k}.   
    \end{split}
\end{equation}
We used $\sum_{1 \leq j_1,..,j_{2k}\leq d}\partial_{x_{j_{2k}}}..\partial_{x_{j_1}}\widehat{u_{j_1...j_{2k}}} = (-i)^{2k}\widehat{u}(x)$ in the last identity. Using inversion formula and integration by parts we further notice that for fixed $j_2,.., j_{2k}$
\begin{align*}
    \partial_{x_k} \widehat{u_{jj_2..j_{2k}}} = \partial_{x_j} \widehat{u_{kj_2..j_{2k}}}.
\end{align*}
This implies that $u_{j_1j_2..j_{2k}} = \partial_{x_{j_1}}\varphi_{j_2..j_{2k}}$ for some smooth function $\varphi_{j_2...j_{2k}}$ with zero average. Furthermore, the periodicity of $\widehat{u_{j_1..j_{2k}}}$ along with its zero average implies that $\varphi_{j_2..j_{2k}}$ is $2\pi$-periodic in each variable. Since, $u_{j_1j_2..j_{2k}}$ is symmetric w.r.t. to $j_1$ and $j_2$ we get $\partial_{x_k}\varphi_{jj_3..j_{2k}} = \partial_{x_j} \varphi_{kj_3..j_{2k}}$. This implies that $\varphi_{j_2..j_{2k}}= \partial_{x_{j_2}}\xi_{j_3..j_{2k}}$ for some smooth $2\pi$-periodic function $\xi_{j_3...j_{2k}}$ whose average is zero. Using this argument iteratively we get a smooth $2\pi$-periodic function $\psi(x)$ with zero average such that
\begin{equation}\label{2.17}
    \widehat{u_{j_1..j_{2k}}}(x) = \partial_{x_{j_{2k}}}...\partial_{x_{j_1}} \psi(x).
\end{equation}
Using \eqref{2.17} in \eqref{2.14} and \eqref{2.16}, we get expressions \eqref{2.11} and \eqref{2.12} respectively.
\end{proof}

\section{Hardy-type inequalities on a torus}\label{sec:hardytorus}
In this section, we will prove Hardy-type inequalities for the operators $\Delta^m$ and $\nabla (\Delta^m)$ for non-negative integers $m$ on the torus $Q_d$. We begin by proving a weighted Hardy inequality for the gradient. 

\begin{theorem}[Weighted Hardy inequality]\label{thm3.1}
Let $\psi \in C^\infty(\overline{Q_d})$ all of whose derivatives are $2\pi$-periodic in each variable. Further assume that $\psi$ has zero average. Let $k$ be a non-positive integer. Then for $d > -2k +2$ we have
\begin{equation}\label{3.1}
    H(k, d) \int_{Q_d} |\psi(x)|^2 \omega(x)^{k-1} dx \leq \int_{Q_d} |\nabla \psi(x)|^2 \omega(x)^{k} dx,
\end{equation}
where $\omega(x) := \sum_j \sin^2(x_j/2)$, 
\begin{equation}\label{3.2}
    H(k, d)^{-1} := \sum_{j=0}^{-k} d^j C_1(k + j, d) \prod_{i=0}^{j-1}C_2(k +i, d) + d^{-k}\prod_{i=0}^{-k}C_2(k+i, d),
\end{equation}
 
\begin{equation}
    C_1(k, d) := 16/(d+2k-2)^2 \hspace{9pt} \text{and} \hspace{9pt} C_2(k, d) := (3d+2k-2)/d(d+2k-2).
\end{equation}
\end{theorem}
\begin{remark}
Note that $C_1(k, d) \sim 1/d^2 , C_2(k,d) \sim 1/d$. This implies that $H(k, d)\sim d$ as $ d \rightarrow \infty$.
\end{remark}

\begin{proof}
The proof goes via expansion of squares. Let $F = (F_1,..,F_d)$ be a smooth real-valued vector field on $Q_d$ which is $2\pi-$ periodic in each variable. Let $\omega_{\epsilon} := \omega+\epsilon^2$ for $\epsilon \neq 0$. Consider the following square, for non-positive integer $k$ and a real parameter $\beta$
\begin{align*}
    |\omega_\epsilon^{k/2} \nabla \psi  + \beta \omega_\epsilon^{k/2} \psi F|^2 &= \omega_\epsilon^{k}|\nabla \psi|^2 + \beta^2 \omega_\epsilon^{k}|F|^2 |\psi|^2 + 2\omega_\epsilon^{k} \beta  \text{Re}(\overline{\psi}\nabla \psi \cdot F) \\
    &= \omega_\epsilon^{k}|\nabla \psi|^2 + \beta^2 \omega_\epsilon^{k}|F|^2 |\psi|^2 + \beta \omega_\epsilon^{k} F \cdot \nabla |\psi|^2.
\end{align*}
Integrating both sides and applying integration by parts, we obtain
\begin{align*}
    0 \leq \int_{Q_d}|\omega_\epsilon^{k/2} \nabla \psi  + \beta \omega_\epsilon^{k/2} \psi F|^2 dx &= \int_{Q_d} \omega_\epsilon^{k} |\nabla \psi|^2 dx \\
    &+ \int_{Q_d} (\beta^2 \omega_\epsilon^{k}|F|^2 - \beta \text{div}(\omega_\epsilon^{k}F))|\psi|^2 dx, 
\end{align*}
which implies 
\begin{equation}\label{3.4}
    \int_{Q_d} |\nabla \psi|^2 \omega_\epsilon^{k} dx \geq \int_{Q_d} \Big(\beta\text{div}(\omega_\epsilon^{k}F)-\beta^2 \omega_\epsilon^{k}|F|^2\Big)|\psi|^2 dx.
\end{equation}
Let $F_j = \frac{\sin x_j}{\omega_\epsilon}$, then 
\begin{align*}
    |F|^2 &= \frac{4}{\omega_\epsilon^2}(\omega-\sum_j \sin^4(x_j/2)) \hspace{5pt}\text{and}\\
    \text{div}(\omega_\epsilon^{k}F)&= d\omega_\epsilon^{k-1} + 2(k-1)\omega\omega_\epsilon^{k-2} -2\omega\omega_\epsilon^{k-1} - 2(k-1)\omega_\epsilon^{k-2}\sum_i \sin^4(x_i/2).
\end{align*}
Using above expressions we obtain
\begin{align*}
    \beta\text{div}(\omega_\epsilon^{k}F)-\beta^2 \omega_\epsilon^{k}|F|^2 &= d\beta\omega_\epsilon^{k-1} + \Big(2\beta(k-1)-4\beta^2\Big)\omega \omega_\epsilon^{k-2}\\
    &- \Big(2\beta(k-1)-4\beta^2\Big)\omega_\epsilon^{k-2}\sum_j \sin^4(x_j/2) -2\beta \omega\omega_\epsilon^{k-1}.
\end{align*}
Plugging the above identity in \eqref{3.4}, and taking limit $\epsilon \rightarrow 0$, we get for $ d>-2k + 2$
\begin{align*}
    \int_{Q_d} |\nabla \psi|^2 \omega^{k} dx &\geq \Big(-4\beta^2 + \beta(d+2k-2)\Big) \int_{Q_d} |\psi|^2 \omega^{k-1} dx\\
    &- \int_{Q_d}\Bigg(\Big(2\beta(k-1) -4\beta^2\Big)\sum_j \sin^4(x_j/2)/\omega^2 + 2\beta \Bigg) \omega^{k}|\psi|^2 dx.
\end{align*}
Choosing $\beta = (d+2k-2)/8$ with the aim of maximizing $-4\beta^2 + \beta(d+2k-2)$, and using the estimate $d\sum_j \sin^4(x_j/2) \geq \omega^2$, we obtain
\begin{align*}
   \int_{Q_d} |\psi|^2 \omega^{k-1} dx &\leq \frac{16}{(d+2k-2)^2} \int_{Q_d}|\nabla \psi|^2 \omega^{k} dx + \frac{3d+2k-2}{d(d+2k-2)}\int_{Q_d}|\psi|^2 \omega^{k} dx\\
   &=: C_1(k, d) \int_{Q_d}|\nabla \psi|^2 \omega^{k} dx + C_2(k, d) \int_{Q_d}|\psi|^2 \omega^{k} dx
\end{align*}
Applying the above inequality inductively w.r.t. $k$ and using $\omega(x) \leq d$, we get 
\begin{align*}
    \int_{Q_d} |\nabla \psi|^2 \omega^{k} dx &\leq \sum_{j=0}^{-k} C_1(k + j, d) \prod_{i=0}^{j-1}C_2(k+i, d) \int_{Q_d} |\nabla \psi|^2 \omega^{k+j} dx\\
    &+ \prod_{i=0}^{-k}C_2(k+i, d)\int_{Q_d} |\psi|^2 dx\\
    & \leq \sum_{j=0}^{-k} d^j C_1(k + j, d) \prod_{i=0}^{j-1}C_2(k+i, d) \int_{Q_d} |\nabla \psi|^2 \omega^{k} dx \\
    &+ \prod_{i=0}^{-k}C_2(k+i, d)\int_{Q_d} |\psi|^2 dx.
\end{align*}
Finally, using the Poincare-Friedrichs inequality: $\int_{Q_d} |\psi|^2 \leq \int_{Q_d} |\nabla \psi|^2$ along with $\omega(x)^k \geq d^k$, we get the desired result.
\end{proof}

\begin{corollary}[Hardy inequality]\label{cor3.3}
Let $\psi \in C^\infty(\overline{Q_d})$ all of whose derivatives are $2\pi$-periodic in each variable. Further assume that $\psi$ has zero average. Then, for $d \geq 3$, we have
\begin{equation}\label{3.5}
    \frac{d(d-2)^2}{3d^2+8d+4} \int_{Q_d} \frac{|\psi(x)|^2}{\sum_j \sin^2(x_j/2)} dx \leq \int_{Q_d} |\nabla \psi(x)|^2 dx.
\end{equation}
\end{corollary}

\begin{proof}
Applying Theorem \ref{thm3.1} for $k=0$ and noting that $H(0, d)=\frac{d(d-2)^2}{3d^2+8d+4} $, the proof is complete.
\end{proof}

In the next lemma, we prove a two-parameter family of inequalities from which we will derive a weighted Rellich and Hardy-Rellich type inequalities on the torus.
\begin{lemma}\label{lem3.4}
Let $\psi \in C^\infty(\overline{Q_d})$ all of whose derivatives are $2\pi$-periodic in each variable. Let $\alpha \leq 0$ and $ d> -4\alpha+4$. Then for real parameters $\beta, \gamma$ satisfying $\beta^2 -\beta(2\alpha-1) \geq 0$ we have
\begin{equation}\label{3.6}
    \begin{split}
        \int_{Q_d} \omega(x)^{2\alpha} |\Delta \psi(x)|^2 dx &\geq (2\gamma -\beta(d+4\beta-4\alpha)) \int_{Q_d} \omega(x)^{2\alpha-1}|\nabla \psi(x)|^2 dx\\
        &+ \frac{\gamma}{2}\Big((2\beta-2\alpha+1)(d+4\alpha-4) -2\gamma \Big) \int_{Q_d} \omega(x)^{2\alpha-2}|\psi(x)|^2 dx\\
        &+ E(x),
    \end{split}
\end{equation}
where $\omega (x) := \sum_i \sin^2(x_i/2)$ and 
\begin{equation}\label{3.7}
    \begin{split}
        E(x) &:= 2\beta \int_{Q_d} \Big(1 +(2\beta -2\alpha+1)\sum_i \sin^4(x_i/2)/\omega^2 \Big)\omega(x)^{2\alpha} |\nabla \psi(x)|^2 dx\\
        &-4\beta \int_{Q_d}\omega(x)^{2\alpha} \sum_i \frac{\sin^2(x_i/2)}{\omega}  |\partial_{x_i}\psi(x)|^2 dx\\
        &-\gamma(2\beta -2\alpha +1)\int_{Q_d} \Big(1 + 2(\alpha-1)\sum_i\sin^4(x_i/2)/\omega^2\Big)\omega(x)^{2\alpha-1}|\psi(x)|^2 dx.
    \end{split}
\end{equation}

\end{lemma}

\begin{proof}
Let $F = (F_1,.., F_d)$ and $f$ be a smooth $2\pi$-periodic real-valued vector field and a scalar function respectively. Let $\omega_\epsilon := \omega+ \epsilon^2$ and $\beta, \gamma \in \R$. We expand the following square
\begin{equation}\label{3.8}
    \begin{split}
        Q(\psi) :=|\omega_\epsilon^\alpha \Delta \psi  + \beta \omega_\epsilon^\alpha F \cdot \nabla \psi + \gamma \omega_\epsilon^\alpha f \psi|^2 &= \omega_\epsilon^{2\alpha}|\Delta \psi|^2 + \beta^2 \omega_\epsilon^{2\alpha}|F \cdot \nabla \psi|^2 + \gamma^2 \omega_\epsilon^{2\alpha}f^2 |\psi|^2\\
        &+ 2\beta \omega_\epsilon^{2\alpha} \text{Re}(F \cdot \nabla \psi \Delta \overline{\psi})+ 2\gamma \omega_\epsilon^{2\alpha} \text{Re}(f \overline{\psi} \Delta \psi) \\
        &+ 2\beta \gamma \omega_\epsilon^{2\alpha} \text{Re}(f F\cdot \overline{\psi} \nabla \psi). 
    \end{split}
\end{equation}
Applying chain rule and integration by parts we obtain
\begin{align*}
    2 \beta \gamma \int_{Q_d}\omega_\epsilon^{2\alpha} \text{Re}(f F\cdot \overline{\psi} \nabla \psi) dx &= -\beta \gamma \int_{Q_d} \text{div}(\omega_\epsilon^{2\alpha}fF)|\psi|^2 dx,\\
    2\gamma \int_{Q_d} \omega_{\epsilon}^{2\alpha} \text{Re}(f \overline{\psi} \Delta \psi) dx &= -2\gamma \int_{Q_d} (\omega_\epsilon^{2\alpha}f)|\nabla \psi|^2 dx -2\gamma \int_{Q_d} \partial_{x_i} (\omega_\epsilon^{2\alpha} f) \text{Re}(\overline{\psi} \partial_{x_i}\psi) dx\\
    &=-2\gamma \int_{Q_d}(\omega_\epsilon^{2\alpha}f)|\nabla \psi|^2 dx + \gamma \int_{Q_d} \Delta(\omega_\epsilon^{2\alpha}f)|\psi|^2 dx,
\end{align*} 
and 
\begin{equation}\label{3.9}
    \begin{split}
        2\beta \int_{Q_d} \omega_\epsilon^{2\alpha} \text{Re}(F \cdot \nabla \psi \Delta \overline{\psi}) dx &= -2\beta \sum_{i, j}\int_{Q_d} (\omega_\epsilon^{2\alpha} F_i) \partial_{x_j}(\partial_{x_i} \psi) \overline{\partial_{x_j}\psi} dx \\ 
        &-2\beta \sum_{i, j} \int_{Q_d} \partial_{x_j}(\omega_\epsilon^{2\alpha}F_i) \partial_{x_i} \psi \partial_{x_j} \overline{\psi} dx\\
        &=\beta \int_{Q_d} \text{div}(\omega_\epsilon^{2\alpha}F) |\nabla \psi|^2 dx \\
        &- 2\beta \sum_{i, j} \int_{Q_d} \partial_{x_j}(\omega_\epsilon^{2\alpha}F_i) \partial_{x_i} \psi \partial_{x_j} \overline{\psi} dx.     
    \end{split}
\end{equation}
Let $F_i(x) := \sin x_i/\omega_\epsilon$ and $f(x) := 1/\omega_\epsilon $. Then we have
\begin{align*}
    \partial_{x_j}(\omega_\epsilon^{2\alpha}F_i) &= \omega_\epsilon^{2\alpha-1}(1-2\sin^2(x_i/2)) \delta_{ij} + (2\alpha-1)/2\omega_\epsilon^{2\alpha-2}\sin x_i \sin x_j,\\
    \partial_{x_i}(\omega_\epsilon^{2\alpha} f F_i) &= \omega_\epsilon^{2\alpha-2}(1-2\sin^2(x_i/2)) + 2(2\alpha-2)\omega_\epsilon^{2\alpha-3}(\sin^2(x_i/2)-\sin^4(x_i/2)),\\
    \partial^2_{x_i^2} (\omega_\epsilon^{2\alpha}f)&= (2\alpha-1)/2\omega_\epsilon^{2\alpha-2}(1-2\sin^2(x_i/2)) \\
    &+ (2\alpha-1)(2\alpha-2)\omega_\epsilon^{2\alpha-3}(\sin^2(x_i/2)-\sin^4(x_i/2)).
\end{align*}
Plugging the above identities in \eqref{3.9}, we get,
\begin{equation}\label{3.10}
    \begin{split}
        &2 \beta \gamma \int_{Q_d}\omega_\epsilon^{2\alpha} \text{Re}(f F\cdot \overline{\psi} \nabla \psi) dx\\ &= -\beta \gamma\int_{Q_d}  \Big(d+ (4\alpha -4)\omega/\omega_\epsilon \Big) \omega_\epsilon^{2\alpha-2}|\psi|^2 dx\\
        &+ 2\beta \gamma \int_{Q_d} \Big(\omega/\omega_\epsilon + 2(\alpha-1)\sum_i \sin^4(x_i/2)/\omega_\epsilon^2\Big)\omega_\epsilon^{2\alpha-1}|\psi|^2 dx,
    \end{split}
\end{equation}

\begin{equation}\label{3.11}
    \begin{split}
        &2\gamma \int_{Q_d} \omega_{\epsilon}^{2\alpha} \text{Re}(f \overline{\psi} \Delta \psi) dx \\
        &= -2\gamma \int_{Q_d}\omega_\epsilon^{2\alpha-1}|\nabla \psi|^2 dx + \gamma(2\alpha-1)\int_{Q_d} \Big(d/2 + 2(\alpha-1)\omega/\omega_\epsilon \Big)\omega_\epsilon^{2\alpha-2}|\psi|^2 dx\\
        &-\gamma(2\alpha-1)\int_{Q_d}\Big(\omega/\omega_\epsilon +  2(\alpha-1)\sum_i \sin^4(x_i/2)/\omega_\epsilon^2\Big)\omega_\epsilon^{2\alpha-1}|\psi|^2 dx,
    \end{split}
\end{equation}
and
\begin{equation}\label{3.12}
    \begin{split}
       &2\beta \int_{Q_d} \omega_\epsilon^{2\alpha} \text{Re}(F \cdot \nabla \psi \Delta \overline{\psi}) dx\\
       &= \int_{Q_d}\Big(d \beta + 2\beta(2\alpha-1)\omega/\omega_\epsilon-2\beta\Big) \omega_\epsilon^{2\alpha-1}|\nabla \psi|^2 dx\\
       &- \beta(2\alpha-1)\int_{Q_d}\omega_\epsilon^{2\alpha} |F \cdot \nabla \psi|^2 dx + 4\beta \int_{Q_d}\omega_\epsilon^{2\alpha} \sum_i \frac{\sin^2(x_i/2)}{\omega_\epsilon}  |\partial_{x_i}\psi|^2 dx \\
       &- 2\beta \int_{Q_d} \Big(\omega/\omega_\epsilon +(2\alpha-1)\sum_i \sin^4(x_i/2)/\omega_\epsilon^2 \Big)\omega_\epsilon^{2\alpha} |\nabla \psi|^2 dx.  
    \end{split}
\end{equation}
Integrating both sides in \eqref{3.8}, and using identities \eqref{3.10}-\eqref{3.12} we obtain
\begin{align*}
    \int_{Q_d}Q(\psi) dx &= \int_{Q_d} \omega_\epsilon^{2\alpha} |\Delta \psi|^2 dx + \int_{Q_d} \Big(-2\gamma + \beta \big(d-2 + (4\alpha-2)\omega/\omega_\epsilon\big)\Big) \omega_{\epsilon}^{2\alpha-1}|\nabla \psi|^2 dx\\
    &+ \int_{Q_d} \frac{\gamma}{2}\Big(2\gamma + (2\alpha-2\beta-1)(d+4(\alpha-1)\omega/\omega_\epsilon)\Big) \omega^{2\alpha-2} |\psi|^2 dx\\
    &+ \Big(\beta^2 -\beta(2\alpha-1)\Big)\int_{Q_d}\omega_\epsilon^{2\alpha} |F \cdot \nabla \psi|^2 dx + E(x), 
\end{align*}
where 
\begin{align*}
    E(x) &:= 4\beta \int_{Q_d}\omega_\epsilon^{2\alpha} \sum_i \frac{\sin^2(x_i/2)}{\omega_\epsilon}  |\partial_{x_i}\psi|^2 dx \\
    &- 2\beta \int_{Q_d} \Big(\omega/\omega_\epsilon +(2\alpha-1)\sum_i \sin^4(x_i/2)/\omega_\epsilon^2 \Big)\omega_\epsilon^{2\alpha} |\nabla \psi|^2 dx\\
    &+ \gamma(2\beta -2\alpha +1)\int_{Q_d} \Big(\omega/\omega_\epsilon + 2(\alpha-1)\sum_i\sin^4(x_i/2)/\omega_\epsilon^2\Big)\omega_\epsilon^{2\alpha-1}|\psi|^2 dx.
\end{align*}
Further using the non-negativity of $\int_{Q_d} Q(\psi) dx$ along with Cauchy's inequality($|F \cdot \nabla \psi| \leq |F||\nabla \psi|$), we get, for $\beta^2 - \beta(2\alpha-1) \geq 0$
\begin{equation}\label{3.13}
    \begin{split}
        \int_{Q_d} \omega_\epsilon^{2\alpha} |\Delta \psi|^2 dx & \geq \int_{Q_d} \Big(2\gamma - \beta (d-2+ 2(2\beta-2\alpha+1) \omega/\omega_\epsilon)\Big) \omega_{\epsilon}^{2\alpha-1}|\nabla \psi|^2 dx\\
        &+ \int_{Q_d} \frac{\gamma}{2}\Big((2\beta - 2\alpha+1)(d+4(\alpha-1)\omega/\omega_\epsilon) -2\gamma\Big) \omega^{2\alpha-2} |\psi|^2 dx\\
        &+ E_1(x),
    \end{split}
\end{equation}
where 
\begin{equation}\label{3.14}
    \begin{split}
        E_1(x) &:= -4\beta \int_{Q_d}\omega_\epsilon^{2\alpha} \sum_i \frac{\sin^2(x_i/2)}{\omega_\epsilon}  |\partial_{x_i}\psi|^2 dx \\
        &+ 2\beta \int_{Q_d} \Big(\omega/\omega_\epsilon +(2\beta -2\alpha+1)\sum_i \sin^4(x_i/2)/\omega_\epsilon^2 \Big)\omega_\epsilon^{2\alpha} |\nabla \psi|^2 dx\\
        &-\gamma(2\beta -2\alpha +1)\int_{Q_d} \Big(\omega/\omega_\epsilon + 2(\alpha-1)\sum_i\sin^4(x_i/2)/\omega_\epsilon^2\Big)\omega_\epsilon^{2\alpha-1}|\psi|^2 dx.    
    \end{split}
\end{equation}
Taking limit $\epsilon \rightarrow 0$ in \eqref{3.13}, we get for $d > -4\alpha+4$
\begin{equation}\label{3.15}
    \begin{split}
        \int_{Q_d} \omega^{2\alpha} |\Delta \psi|^2 dx &\geq (2\gamma -\beta(d+4\beta-4\alpha)) \int_{Q_d} \omega^{2\alpha-1}|\nabla \psi|^2 dx\\
        &+ \frac{\gamma}{2}\Big((2\beta-2\alpha+1)(d+4\alpha-4) -2\gamma \Big) \int_{Q_d} \omega^{2\alpha-2}|\psi|^2 dx + E_2(x),
    \end{split}
\end{equation}
where
\begin{equation}\label{3.16}
    \begin{split}
        E_2(x) &:= -4\beta \int_{Q_d}\omega^{2\alpha} \sum_i \frac{\sin^2(x_i/2)}{\omega}  |\partial_{x_i}\psi|^2 dx\\
        &+ 2\beta \int_{Q_d} \Big(1 +(2\beta -2\alpha+1)\sum_i \sin^4(x_i/2)/\omega^2 \Big)\omega^{2\alpha} |\nabla \psi|^2 dx\\
        &-\gamma(2\beta -2\alpha +1)\int_{Q_d} \Big(1 + 2(\alpha-1)\sum_i\sin^4(x_i/2)/\omega^2\Big)\omega^{2\alpha-1}|\psi|^2 dx.
    \end{split}
\end{equation}

\end{proof}

Next, we prove a weighted Hardy-Rellich type inequality on the torus $Q_d$.
\begin{theorem}[Weighted Hardy-Rellich inequality]\label{thm3.5}
Let $\psi \in C^\infty(\overline{Q_d})$ all of whose derivatives are $2\pi$-periodic in each variable. Further assume that $\psi$ has zero average. Let $k$ be a non-positive integer. Then for $d \geq -6k +8$ we have
\begin{equation}\label{3.17}
    HR(k, d) \int_{Q_d} |\nabla \psi(x)|^2 \omega(x)^{k-1} dx \leq \int_{Q_d} |\Delta \psi(x)|^2 \omega(x)^k dx,
\end{equation}
where $\omega(x) := \sum_j \sin^2(x_j/2)$, 
\begin{equation}\label{3.18}
    HR(k, d)^{-1} := \sum_{j=0}^{-k} d^j C_1(k + j, d) \prod_{i=0}^{j-1}C_2(k + i, d) + d^{-k}\prod_{i=0}^{-k}C_2(k+i, d),
\end{equation}

\begin{equation}
    C_1(k, d) := 16/(d-2k)^2 \hspace{9pt} \text{and} \hspace{9pt} C_2(k, d) := (3d-2k+4)/d(d-2k).
\end{equation}
\end{theorem}

\begin{remark}
Note that $C_1(k, d) \sim 1/d^2$ and $C_2(k, d) \sim 1/d$, which implies that $HR(k, d) \sim d$ as $d \rightarrow \infty$.
\end{remark}

\begin{proof}
We begin by writing down inequality \eqref{3.6}, with $\alpha \leq 0$, $\beta, \gamma$ being real numbers satisfying $\beta^2 -\beta(2\alpha-1) \geq 0$ and $d >-4\alpha+4$
\begin{equation}\label{3.20}
    \begin{split}
        \int_{Q_d} \omega^{2\alpha} |\Delta \psi|^2 dx &\geq (2\gamma -\beta(d+4\beta-4\alpha)) \int_{Q_d} \omega^{2\alpha-1}|\nabla \psi|^2 dx\\
        &+ \frac{\gamma}{2}\Big((2\beta-2\alpha+1)(d+4\alpha-4) -2\gamma \Big) \int_{Q_d} \omega^{2\alpha-2}|\psi|^2 dx + E(x).
    \end{split}
\end{equation}
Next, we choose $\gamma =0$ and $\beta = -(d-4\alpha)/8$. 
Note that for this choice of $\gamma$, the second term on the RHS of \eqref{3.20} vanishes and $\beta = -(d-4\alpha)/8 $ maximizes the coefficient of the first term on the RHS of \eqref{3.20} after taking $\gamma =0$. Also note that condition $\beta^2-\beta(2\alpha-1) \geq 0$ gives an additional constraint $d \geq -12\alpha+8$. For this choice of parameters \eqref{3.20} becomes
\begin{equation}\label{3.21}
    \int_{Q_d} \omega^{2\alpha}|\Delta \psi|^2 dx \geq \frac{(d-4\alpha)^2}{16} \int_{Q_d} \omega^{2\alpha-1} |\nabla \psi|^2 dx + E(x),
\end{equation}
and $E(x)$ becomes
\begin{align*}
    E(x) &= 2\beta \int_{Q_d} \Big(1 +(2\beta -2\alpha+1)\sum_i \sin^4(x_i/2)/\omega^2 \Big)\omega^{2\alpha} |\nabla \psi|^2 dx\\
    &-4\beta \int_{Q_d}\omega^{2\alpha} \sum_i \frac{\sin^2(x_i/2)}{\omega}  |\partial_{x_i}\psi|^2 dx\\
    &= -\frac{(d-4\alpha)}{4} \int_{Q_d} \omega^{2\alpha} |\nabla \psi|^2 dx + \frac{(d-4\alpha)(d+4\alpha-4)}{16} \int_{Q_d} \frac{\sum_i \sin^4(x_i/2)}{\omega^2} \omega^{2\alpha} |\nabla \psi|^2 dx\\
    &+(d-4\alpha)/2 \int_{Q_d}\omega^{2\alpha} \sum_i \frac{\sin^2(x_i/2)}{\omega}  |\partial_{x_i}\psi|^2 dx\\ 
    & \geq -\frac{(d-4\alpha)(3d-4\alpha +4)}{16 d}\int_{Q_d} \omega^{2\alpha}|\nabla \psi|^2 dx\\
\end{align*}
In the last inequality we used $d \sum_i \sin^4(x_i/2) \geq \omega^2$ and $\alpha \leq 0$ to bound the last integral by zero from below. Using this lower bound on $E(x)$ in \eqref{3.21} and taking $\alpha := k/2$ we obtain
\begin{equation}\label{3.22}
    \int_{Q_d} |\nabla \psi|^2 \omega^{k-1} dx \leq  \frac{16}{(d-2k)^2}\int_{Q_d} |\Delta \psi|^2 \omega^{k} dx + \frac{3d-2k+4}{d(d-2k)}\int_{Q_d}  |\nabla \psi|^2 \omega^{k} dx.
\end{equation}
Applying inequality \eqref{3.22} inductively w.r.t $k$ and using the bound $w(x) \leq d$ we get

\begin{align*}
    \int_{Q_d} |\nabla \psi|^2 \omega^{k-1} dx &\leq \sum_{j=0}^{-k} d^j C_1(k + j, d) \prod_{i=0}^{j-1}C_2(k+i, d) \int_{Q_d} |\Delta \psi|^2 \omega^k dx \\
    &+ \prod_{i=0}^{-k}C_2(k+i, d)\int_{Q_d} |\nabla \psi|^2 dx,
\end{align*}
where $C_1(k, d) := 16/(d-2k)^2$ and $C_2(k, d) := (3d-2k+4)/d(d-2k)$.\\
Using $\omega^{k} \geq d^{k} $ and $\int_{Q_d}|\nabla \psi|^2 dx \leq \int_{Q_d}|\Delta \psi|^2 dx$ in the above inequality completes the proof.
\end{proof}

\begin{corollary}[Hardy-Rellich inequality]\label{cor3.7}
Let $\psi \in C^\infty(\overline{Q_d})$ all of whose derivatives are $2\pi$-periodic in each variable. Further assume that $\psi$ has zero average. Then for $d \geq 8$, we have
\begin{equation}\label{3.23}
    \frac{d^2}{3d+20} \int_{Q_d} \frac{|\nabla \psi(x)|^2}{\sum_j \sin^2(x_j/2)} dx \leq \int_{Q_d} |\Delta \psi(x)|^2 dx.
\end{equation}
\end{corollary}

\begin{proof}
Inequality \eqref{3.23} follows directly from Theorem \ref{thm3.5} by taking $k=0$ and observing that $HR(0, d) = d^2/(3d+20)$.
\end{proof}

Next, we derive a Rellich inequality on the torus $Q_d$ from Lemma \ref{lem3.4}.
\begin{theorem}[Weighted Rellich inequality]\label{thm3.8}
Let $\psi \in C^\infty(\overline{Q_d})$ all of whose derivatives are $2\pi$-periodic in each variable. Further assume that $\psi$ has zero average. Let $k$ be a non-positive integer. Then, for $d > -2k +4$, we have
\begin{equation}\label{3.24}
    R(k, d) \int_{Q_d} |\psi(x)|^2 \omega(x)^{k-2} dx \leq \int_{Q_d} |\Delta \psi(x)|^2 \omega(x)^k dx,
\end{equation}
where $\omega(x) := \sum_j \sin^2(x_j/2)$, 
\begin{equation}\label{3.25}
    R(k, d) :=  \frac{(d-2k)^2(d+2k-4)^2}{256\Bigg(1+HR(k, d)^{-1}\Big(d C_1(k/2, d)+d C_2(k/2, d)H(k, d)^{-1}\Big)\Bigg)},
\end{equation}
non-negative constants $C_1(\alpha, d), C_2(\alpha, d)$ are given by
\begin{align*}
    C_1(\alpha, d) &:= \frac{2\beta(d-2\beta + 2\alpha-1)}{d},\\ C_2(\alpha, d) &:= \frac{\beta(d+4\beta-4\alpha)(d+2\alpha-2)(2\beta-2\alpha+1)}{2d},
\end{align*}
and $\beta := 1/8(-4 + 8\alpha + \sqrt{2}\sqrt{d^2 -4d + 16\alpha^2 -16 \alpha + 8}) \geq 0$. 
\end{theorem}

\begin{remark}
Note that $R(k, d) \sim d^2$ as $d \rightarrow \infty$, since $H(k, d)\sim d, HR(k, d) \sim d$ and $C_1(k/2, d) \sim d$, $C_2(k/2, d)\sim d^3$. 
\end{remark}
\begin{proof}
Let $\alpha \leq 0$ and $\beta, \gamma$ be real numbers satisfying $\beta^2-\beta(2\alpha-1) \geq 0$. Then we have inequality \eqref{3.6}, that is
\begin{equation}\label{3.26}
    \begin{split}
        \int_{Q_d} \omega^{2\alpha} |\Delta \psi|^2 dx &\geq (2\gamma -\beta(d+4\beta-4\alpha)) \int_{Q_d} \omega^{2\alpha-1}|\nabla \psi|^2 dx\\
        &+ \frac{\gamma}{2}\Big((2\beta-2\alpha+1)(d+4\alpha-4) -2\gamma \Big) \int_{Q_d} \omega^{2\alpha-2}|\psi|^2 dx + E(x),
    \end{split}
\end{equation}
where $E(x)$ is as defined in \eqref{3.7}.  We first choose $2\gamma = \beta(d+4\beta-4\alpha)$ for which the first term in RHS of \eqref{3.26} vanishes. We further choose $\beta = 1/8(-4 + 8\alpha + \sqrt{2}\sqrt{d^2 -4d + 16\alpha^2 -16 \alpha + 8})$ with the aim of maximizing the coefficient of the second term in the RHS of \eqref{3.26}. This choice of parameters implies that
\begin{equation}\label{3.27}
   \int_{Q_d} |\Delta \psi|^2 \omega^{2\alpha}  dx  \geq \frac{(d-4\alpha)^2(d+4\alpha-4)^2}{256}\int_{Q_d}  |\psi|^2 \omega^{2\alpha-2}dx + E(x).
\end{equation}
Note that $\beta, \gamma \geq 0$ and $\beta^2-\beta(2\alpha-1)\geq 0$ under the condition $d > -4\alpha+4$. Using the estimates $\omega \geq \sin^2(x_i/2)$ and $d \sum_i \sin^4(x_i/2) \geq \omega^2$, we get
\begin{align*}
    E(x) \geq &-\frac{2\beta(d -2\beta+2\alpha-1)}{d}\int_{Q_d}\omega^{2\alpha}|\nabla \psi|^2 dx\\
    &- \frac{\gamma(d+ 2\alpha-2)(2\beta -2\alpha +1)}{d} \int_{Q_d}\omega^{2\alpha-1}|\psi|^2 dx.
\end{align*}
Above inequality along with \eqref{3.27} gives
\begin{equation}
    \begin{split}
        \frac{(d-4\alpha)^2(d+4\alpha-4)^2}{256} \int_{Q_d} |\psi|^2 \omega^{2\alpha-2} dx \leq \int_{Q_d} |\Delta \psi|^2 \omega^{2\alpha} dx &+ C_1(\alpha, d) \int_{Q_d} |\nabla \psi|^2 \omega^{2\alpha} dx \\
        &+ C_2(\alpha, d) \int_{Q_d} |\psi|^2 \omega^{2\alpha-1} dx,
    \end{split}
\end{equation}
where $C_1(\alpha, d) := 2\beta(d-2\beta+ 2\alpha-1)/d$ and $C_2(\alpha, d) := \gamma(d+2\alpha-2)(2\beta-2\alpha+1)/d$. Note that $d > -4\alpha+4$ implies $C_1(\alpha, d), C_2(\alpha, d) \geq 0$. From now on, we assume that $2\alpha \in \Z$. Using weighted Hardy inequality \eqref{3.1} and weighted Hardy-Rellich inequality \eqref{3.17} we obtain
\begin{equation}
    \begin{split}
    (d-4\alpha)^2(d+4\alpha-4)^2/256 \int_{Q_d} |\psi|^2 \omega^{2\alpha-2} dx \leq C(\alpha, d)\int_{Q_d} |\Delta \psi|^2 \omega^{2\alpha} dx,
    \end{split}
\end{equation}
where $C(\alpha, d) := 1+HR(2\alpha, d)^{-1}\Big(dC_1(\alpha, d)+dC_2(\alpha, d)H(2\alpha, d)^{-1}\Big)$. \\
Taking $\alpha := k/2$ completes the proof.
\end{proof}

\begin{remark}
Note that Theorem \ref{thm3.8} with $k=0$ gives the Rellich inequality on torus $Q_d$, although unlike the Hardy and Hardy-Rellich inequalities \eqref{3.5}, \eqref{3.23}, the constant in the Rellich inequality has a messy expression.  
\end{remark}

In the next theorem we apply weighted Hardy and Rellich inequality \eqref{3.1}, \eqref{3.24} respectively to prove a Hardy type inequality for the operators $\Delta^m$ and $\nabla( \Delta^m)$.
\begin{theorem}
Let $\psi \in C^\infty(\overline{Q_d})$ all of whose derivatives are $2\pi$-periodic in each variable. Further assume that $\psi$ has zero average. Let $k$ be a non-positive integer and $m$ be a non-negative integer.
\begin{enumerate}
    \item If $ d > -2k+4m$, then
    \begin{equation}\label{3.30}
        C(m, k, d)\int_{Q_d} |\psi(x)|^2 \omega(x)^{k-2m} dx \leq \int_{Q_d} |\Delta^m \psi(x)|^2 \omega(x)^{k} dx,
    \end{equation}
    where 
    \begin{equation}\label{3.31}
        C(m, k, d) := \prod_{i=0}^{m-1} R(k-2i, d),
    \end{equation}
    and $R(k, d)$ is the constant in the Rellich inequality \eqref{3.25}.
    \item If $d > -2k+4m+2$, then
    \begin{equation}\label{3.32}
        \widetilde{C}(m, k, d)\int_{Q_d} |\psi(x)|^2 \omega(x)^{k-(2m+1)} dx  \leq \int_{Q_d} |\nabla(\Delta^m \psi)(x)|^2 \omega(x)^{k} dx,
    \end{equation}
    where 
    \begin{equation}\label{3.33}
        \widetilde{C}(m, k, d) := H(k, d)\prod_{i=0}^{m-1}R(k-2i-1, d),
    \end{equation}
    and $H(k, d)$ is the constant in the Hardy inequality \eqref{3.2}.
\end{enumerate}

\end{theorem}

\begin{proof}
Applying the weighted Rellich inequality \eqref{3.24} with $\psi$ replaced by $\Delta^{m-1} \psi$, we get
\begin{equation}\label{3.34}
    \int_{Q_d} |\Delta^m \psi(x)|^2 \omega(x)^k dx \geq R(k, d) \int_{Q_d} |\Delta^{m-1}\psi(x)|^2 \omega(x)^{k-2} dx.
\end{equation}
Applying inequality \eqref{3.34} inductively w.r.t. $m$ we get our desired result.\\\\
Applying weighted Hardy inequality \eqref{3.1} with $\psi$ replaced by $\Delta^m \psi$ we obtain
\begin{equation}\label{3.35}
    \int_{Q_d} |\nabla (\Delta^m \psi)|^2 \omega^{k-1} \geq H(k, d)\int_{Q_d} |\Delta^m \psi|^2 \omega^{k-1} dx.
\end{equation}
Inequality \eqref{3.35} along with \eqref{3.30} gives inequality \eqref{3.32}.

\end{proof}

\begin{corollary}
Let $\psi \in C^\infty(\overline{Q_d})$ all of whose derivatives are $2\pi$-periodic in each variable. Further assume that $\psi$ has zero average. Let $m$ be a non-negative integer.
\begin{enumerate}
    \item If $d > 4m$, then 
    \begin{equation}\label{3.36}
        C(m, d) \int_{Q_d} \frac{|\psi(x)|^2}{\big(\sum_i \sin^2(x_i/2)\big)^{2m}} dx \leq \int_{Q_d} |\Delta^m \psi(x)|^2 dx,
    \end{equation}
    where $C(m, d) := \prod\limits_{i=0}^{m-1}R(-2i, d)$ and $R(k, d)$ is the constant in the Rellich inequality \eqref{3.25}.
    \item If $d >4m +2$, then
    \begin{equation}\label{3.37}
        \widetilde{C}(m, d) \int_{Q_d} \frac{|\psi(x)|^2}{\big(\sum_i \sin^2(x_i/2)\big)^{2m+1}} dx \leq \int_{Q_d} |\nabla(\Delta^m \psi)(x)|^2 dx,
    \end{equation}
    where $\widetilde{C}(m,d) := H(0, d)\prod\limits_{i=0}^{m-1} R(-2i-1, d)$ and $H(0,d)$ is the constant in the Hardy inequality \eqref{3.2}.
\end{enumerate}

\end{corollary}
\begin{remark}\label{rem3.13}
Note that $C(m, d) \sim d^{2m}$ and $\widetilde{C}(m, d) \sim d^{2m+1}$ as $d \rightarrow \infty$, since $H(k, d) \sim d$ and $R(k, d)\sim d^2$ as $d \rightarrow \infty$.
\end{remark}

\begin{proof}
Putting $k=0$ in \eqref{3.30} and \eqref{3.32} gives inequalities \eqref{3.36} and \eqref{3.37} respectively.
\end{proof}

\section{Proof of the Main results}\label{sec:mainproofs}
\begin{proof}[Proof of theorem \ref{thm1.1}]
From Lemma \ref{lem2.2} we know that there exists a smooth function $\psi$ all of whose derivatives are $2\pi$-periodic and has zero average, such that
\begin{equation}\label{4.1}
    \sum_{n \in \Z^d} \frac{|u(n)|^2}{|n|^{4k+2}} = \int_{Q_d} |\nabla(\Delta^k \psi)|^2 dx,
\end{equation}
and 
\begin{equation}\label{4.2}
    \sum_{n \in \Z^d} |D \Delta^k u(n)|^2 = 4^{2k+1}\int_{Q_d} |\Delta^{2k+1} \psi|^2 \omega(x)^{2k+1}dx,
\end{equation}
for $\omega(x) := \sum_i \sin^2(x_i/2)$.\\
Integration by parts, along with H\"older's inequality, gives
\begin{equation}\label{4.3}
    \begin{split}
        \int_{Q_d} |\nabla(\Delta^k \psi)(x)|^2 dx &= -\int_{Q_d} \Delta^{2k+1} \psi(x)  \psi(x) dx \\
        &\leq \Bigg(\int_{Q_d}|\Delta^{2k+1} \psi(x)|^2 \omega^{2k+1} dx\Bigg)^{1/2} \Bigg(\int_{Q_d} \frac{|\psi|^2}{\omega^{2k+1}}dx\Bigg)^{1/2}.    
    \end{split}
\end{equation}
Inequality \eqref{4.3} along with Hardy-type inequality \eqref{3.37} gives
\begin{equation}\label{4.4}
    \widetilde{C}(k,d) \int_{Q_d} |\nabla(\Delta^k \psi)(x)|^2 dx \leq \int_{Q_d} |\Delta^{2k+1} \psi(x)|^2 \omega^{2k+1} dx.
\end{equation}
where $\widetilde{C}(k, d)$ is as defined by \eqref{3.37}.\\
Inequalities \eqref{4.1} and \eqref{4.2} along with \eqref{4.4} yield
\begin{equation}\label{4.5}
  4^{2k+1} \widetilde{C}(k, d)\sum_{n \in \Z^d} \frac{|u(n)|^2}{|n|^{4k+2}} \leq \sum_{n \in \Z^d} |D\Delta^k u(n)|^2   
\end{equation}
This proves that $C_1(k, d) \geq 4^{2k+1}\widetilde{C}(k, d)$, where $C_1(k, d)$ is the sharp constant in \eqref{1.3}. Next, we will prove an upper bound on $C_1(k, d)$. Let $v(n) := \Delta^{k-1} u(n)$ for $u\in C_c(\Z^d)$. Consider
\begin{align*}
    \sum_{n \in \Z^d} |\Delta^k u(n)|^2 = \sum_{n \in \Z^d} |\Delta v(n)|^2 &= \sum_{n \in \Z^d} \Big(\sum_{j=1}^d(2v(n)-v(n-e_j)-v(n+e_j) \Big)^2\\
    &\leq \sum_{n \in \Z^d} \sum_{j, k} \varphi_j(n) \varphi_k(n)\\
    &\leq \sum_{j, k} \Big(\sum_{n \in \Z^d} |\varphi_j(n)|^2\Big)^{1/2}\Big(\sum_{n \in \Z^d} |\varphi_k(n)|^2\Big)^{1/2},
\end{align*}
where $\varphi_j(n) := |2v(n) - v(n-e_j) -v(n+e_j)|$. Applying H\"older's inequality and invariance of $\Z^d$ w.r.t. translations along co-ordinate directions we get, $\sum \limits_{n \in \Z^d}|\varphi_j(n)|^2  \leq 16 \sum\limits_{n \in \Z^d}|v(n)|^2$. Therefore, we have
\begin{align*}
    \sum_{n \in \Z^d} |\Delta^k u(n)|^2 \leq 16d^2 \sum_{n \in \Z^d}|v(n)|^2 = 16 d^2\sum_{n \in \Z^d} |\Delta^{k-1} u(n)|^2.
\end{align*}
Applying the above inequality iteratively, we obtain
\begin{equation}\label{4.6}
    \sum_{n \in \Z^d} |\Delta^k u(n)|^2 \leq 4^{2k} d^{2k} \sum_{n \in \Z^d}|u(n)|^2.  
\end{equation}
Consider a function $u$ defined on $\Z^d$ as follows: $u(n) := 1$ if $|n| = 1$ and $u(n) := 0$ everywhere else. Applying inequality \eqref{4.6}, we get
\begin{align*}
    \sum_{n \in \Z^d} |D (\Delta^k u)(n)|^2 \leq 4d \sum_{n \in \Z^d} |\Delta^k u(n)|^2 &\leq 4^{2k+1} d^{2k+1} \sum_{n \in \Z^d}|u(n)|^2\\
    &= 4^{2k+1} d^{2k+1} \sum_{n \in \Z^d} \frac{|u(n)|^2}{|n|^{4k+2}}. 
\end{align*}
This proves that $C_1(k, d) \leq 4^{2k+1} d^{2k+1}$. Therefore we have $4^{2k+1} \widetilde{C}(k,d) \leq C_1(k, d) \leq 4^{2k+1} d^{2k+1}$. This proves that $C_1(k, d) \sim d^{2k+1}$ as $ d \rightarrow \infty$, since $\widetilde{C}(k, d) \sim d^{2k+1}$ as $d \rightarrow \infty$(see Remark \ref{rem3.13}).
\end{proof}

\begin{proof}[Proof of theorem \ref{thm1.2}]
From Lemma \ref{lem2.3} we get the following identities
\begin{equation}\label{4.7}
    \sum_{n \in \Z^d} \frac{|u(n)|^2}{|n|^{4k}} = \int_{Q_d} |\Delta^k \psi(x)|^2 dx,
\end{equation}

\begin{equation}\label{4.8}
    \sum_{n \in \Z^d} |\Delta^k u(n)|^2 = 4^{2k}\int_{Q_d} |\Delta^k \psi(x)|^2 \omega(x)^{2k} dx. 
\end{equation}
Integration by parts and H\"older's inequality gives
\begin{equation}\label{4.9}
    \begin{split}
        \int_{Q_d} |\Delta^k \psi(x)|^2 dx &= \int_{Q_d} \Delta^{2k}\psi(x) \psi(x) dx \\
        &\leq \Bigg(\int_{Q_d} |\Delta^k \psi(x)|^2 \omega(x)^{2k} dx\Bigg)^{1/2}\Bigg(\int_{Q_d}\frac{|\psi(x)|^2}{\omega(x)^{2k}} dx\Bigg)^{1/2}.    
    \end{split}
\end{equation}
Inequality \eqref{4.9} along with Hardy inequality \eqref{3.36} gives
\begin{equation}\label{4.10}
    C(k, d)\int_{Q_d} |\Delta^k \psi(x)|^2 dx \leq \int_{Q_d} |\Delta^k \psi(x)|^2 \omega(x)^{2k} dx,
\end{equation}
where $C(k, d)$  is defined in \eqref{3.36}. \\
Inequalities \eqref{4.7}, \eqref{4.8} and \eqref{4.10} gives
\begin{equation}
    4^{2k} C(k, d)\sum_{n \in \Z^d} \frac{|u(n)|^2}{|n|^{4k}} \leq \sum_{n \in \Z^d} |\Delta^k u(n)|^2.
\end{equation}
Therefore we have $C_2(k, d) \geq 4^{2k} C(k, d)$, where $C_2(k, d)$ is the sharp constant in the inequality \eqref{1.4}. Consider a function $u(n) :=1 $ if $|n|=1$ and $u(n) := 0$ otherwise. Using inequality \eqref{4.6}, we get
\begin{align*}
    \sum_{n \in \Z^d} |\Delta^k u(n)|^2 \leq 4^{2k}d^{2k} \sum_{n \in \Z^d}\frac{|u(n)|^2}{|n|^{4k}}. 
\end{align*}
This shows that $C_2(k, d) \leq 4^{2k} d^{2k}$. Therefore we have $4^{2k} C(k, d) \leq C_2(k, d) \leq 4^{2k} d^{2k}$, which implies that $C_2(k, d) \sim d^{2k}$ as $d \rightarrow \infty$, since $C(k, d) \sim d^{2k}$ as $ d \rightarrow \infty$(see Remark \ref{rem3.13}).

\end{proof}

\textbf{Acknowledgments.}
We would like to thank Professor Ari Laptev for suggesting the problem and for various valuable discussions on the topic. We would also like to thank Ashvni Narayanan for careful proof reading. Finally, we thank reviewers for their thorough reading and many helpful suggestions. The author is supported by President's Ph.D. Scholarship, Imperial College London.  \\

%\textbf{Data availability:} Data sharing not applicable to this article as no datasets were generated or analysed during the current study.

%%%%%%%%%%%%%%%%%%%%%%%%%%%%%%%%%%%%%%%%%%%
%%%%%%%%%%%%%%%%%%%%%%%%%%%%%%%%%%%%%%%%%%%

\bibliographystyle{amsalpha}

\end{document}